\documentclass[aos,MSNbibl,citesort,dvips]{arximspdf}
\usepackage{dcolumn}
\usepackage{mathrsfs}
\usepackage{graphicx}

\doi{10.1214/10-AOS857}
\volume{39}
\issue{2}
\pubyear{2011}
\firstpage{838}
\lastpage{864}

\makeatletter
\newcolumntype{d}[1]{D{.}{.}{#1}}
\newcommand{\eqref}[1]{(\ref{#1})}
\newcommand{\implies}{\Longrightarrow}
\newtheorem{theorem}{Theorem}[section]
\newproclaim{defn}[theorem]{Definition}
\newtheorem{prop}[theorem]{Proposition}
\newproclaim{rmk}[theorem]{Remark}
\newtheorem{lem}[theorem]{Lemma}
\newproclaim{exmp}[theorem]{Example}
\makeatother

\begin{document}
\begin{frontmatter}

\title{Bayesian analysis of variable-order, reversible~Markov chains\thanksref{T1}}
\runtitle{Variable-order, reversible Markov chains}

\begin{aug}
\author[A]{\fnms{Sergio} \snm{Bacallado}\corref{}\ead[label=e1]{sergiob@stanford.edu}}
\runauthor{S. Bacallado}
\affiliation{Stanford University}
\address[A]{Department of Structural Biology\\
Stanford University\\
Clark Center, S296 \\
Stanford, California 94305\\
USA\\
\printead{e1}} 
\end{aug}
\thankstext{T1}{Supported by a Smith Stanford Graduate
Fellowship and Grants
NIH-R01-GM062868 and NSF DMS-09-00700.}

\received{\smonth{5} \syear{2010}}
\revised{\smonth{9} \syear{2010}}

%
\begin{abstract}
We define a conjugate prior for the reversible Markov chain of order $r$.
The prior arises from a partially exchangeable reinforced random walk,
in the same way that the Beta distribution arises from the exchangeable
Poly\'{a} urn. An extension to variable-order Markov chains is also
derived. We show the utility of this prior in testing the order and
estimating the parameters of a reversible Markov model.
\end{abstract}

%

\end{frontmatter}
%

\section{Introduction}\label{sec1}

Reversible Markov chains are central to a number of fields. They
underlie problems in applied probability like card-shuffling and
queueing networks \cite{AldousFill,Kelly1979fv} and pervade
computational statistics through the many variants of Markov chain
Monte Carlo; in physics, they are natural stochastic models for
time-reversible dynamics. However, the notion of reversibility in
stochastic proscesses with memory is not as widely discussed, and
statistical problems like testing the order of a reversible process
remain a challenge.

We define a conjugate prior for higher-order, reversible Markov chains,
which extends a prior for reversible Markov chains by Diaconis and
Rolles~\cite{Diaconis2006gd}. We begin by defining reversibility in a
more general setting and motivating the significance of higher-order
processes. In Section~\ref{sec2}, we present two graphical representations for
an order-$r$, reversible Markov chain, which are used in Section~\ref{sec3} to
derive the conjugate prior via a random walk with reinforcement. We
dedicate Section~\ref{sec4} to variable-order Markov chains, a family of models
that avoids the curse of dimensionality associated with higher-order
Markov chains, proving essential in certain applications. Finally in
Section~\ref{sec5}, we discuss properties of the prior pertaining to Bayesian
analysis. In examples, we test the extent of memory of a lumped Markov
chain and discretized molecular dynamics trajectories, and compare the
posterior inferences of different models.

\begin{defn}
\label{defReversibility}
A stochastic process $X = X_n, n\in\mathbb{N}$, with distribution $P$ is
called \textit{reversible}, if for any $m> n>0$,
\[
P(X_1,X_2,\ldots ,X_{n}) = P(X_{m-1},X_{m-2},\ldots ,X_{m-n}).
\]
\end{defn}

It is not difficult to show that reversibility implies stationarity
\cite{Kelly1979fv}; if stationarity is given, the above condition need
only be checked for $m=n+1$. Now suppose $X$ is an order-$r$,
irreducible Markov chain taking values in a~finite set $\mathcal{X}$.
We will also apply the term \textit{reversible} to this process when the
stationary chain satisfies the reversibility condition.
\begin{prop}
\label{ReversibilityCondition}
Let $P$ be the stationary law of the order-$r$ Markov chain $X$. If
$P(X_1,\ldots ,X_{r+1}) = P(X_{r+1},\ldots ,X_1)$, then the Markov chain
is reversible.
\end{prop}

\begin{pf}
It is not difficult to check that the hypothesis together with
stationarity imply $P(X_1,\ldots ,X_n)=P(X_n,\ldots ,X_1)$ for any
$n<r+1$. For any $n>r+1$:
\begin{eqnarray*}
P(X_1,\ldots ,X_n) & =& P(X_1,\ldots ,X_{r+1}) \prod_{i=r+2}^n P(X_i |
X_{i-r},\ldots ,X_{i-1}) \\
& =& P(X_1,\ldots ,X_{r+1}) \frac{P(X_2,\ldots ,X_{r+2})}{P(X_2,\ldots
,X_{r+1})}\cdots\frac{P(X_{n-r},\ldots ,X_n)}{P(X_{n-r},\ldots
,X_{n-1})} \\
& =& P(X_n,\ldots ,X_{n-r}) \frac{P(X_{n-1},\ldots
,X_{n-r-1})}{P(X_{n-1},\ldots ,X_{n-r})}\cdots\frac{P(X_{r+1},\ldots
,X_1)}{P(X_{r+1},\ldots ,X_{2})} \\
& =& P(X_n,\ldots ,X_1),
\end{eqnarray*}
where we have used the Markov property, stationarity  and the
hypothesis.
\end{pf}

As a first remark, note that $X_n, n\in\mathbb{N}$, can be represented
as a first-order Markov chain $V_n,n\in\mathbb{N}$, taking values in the
space of sequences $\mathcal{X}^r$. However, the reversibility of $X$
does not imply the reversibility of its first-order representation;
therefore, the analysis of higher-order reversible Markov chains
requires novel techniques. In the following sections, we often use the
first-order representation $V_n,n\in\mathbb{N}$, referring to it
nonetheless as an order-$r$ Markov chain and using the notion of
reversibility associated with the \mbox{order-$r$} Markov chain.

Secondly, we recall that Kolmogorov's criterion is another necessary
and sufficient condition for the reversibility of a Markov chain, which
only depends on the conditional transition probabilities \cite
{Kelly1979fv}. Its equivalence to Definition \ref{defReversibility} in
the higher-order case is proven in the \hyperref[appm]{Appendix}. Kolmogorov's criterion
requires that the probability of traversing any cycle in either
direction is the same. Accordingly, a reversible Markov chain can be
interpreted as a process with no net circulation in space.

%
\begin{figure}
\includegraphics{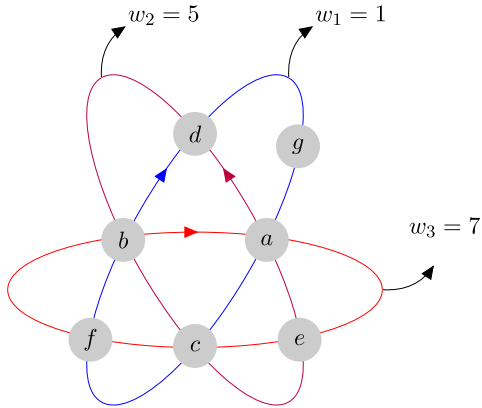}

\caption{A set of weighted circuits on the set $\mathcal{X}=\{
a,b,c,d,e,f,g\}$. In a circuit process started at $u$ in $\mathcal
{X}^r$, we transition on some circuit that contains $u$ with
probability proportional to its weight.}
\label{circuit process picture}

\end{figure}

Reversibility is preserved under certain transformations. For example,
let $X_n, n\in\mathbb{N}$, be a stationary, reversible Markov chain and
consider a finitely valued function, $f(X_n),n\in\mathbb{N}$. It is
easy to check that this process is stationary and reversible, even
though it may not be a Markov chain of any finite order. Functions or
projections of reversible Markov chains appear under different guises
in physics and other fields, and in many cases the effects of memory
subside with time, motivating the use of finite order models. The
problems of determining the order and estimating the parameters of
Markov models have been studied extensively; here, we address these
problems with the constraint of reversibility.

\section{Graphical representations of reversible Markov chains}\label{sec2}

For any sequence $u\in\mathcal{X}^s$, let $u^*$ be its inverse,
$\mathrm{A}(u)$ the subsequence obtained by deleting its last element
and $\Omega(u)$ the one obtained by deleting its first element. We call
$u_1,u_2,\ldots ,u_n$ with $u_i\in\mathcal{X}^s$ an \textit{admissible
path} if $\Omega(u_i)=\mathrm{A}(u_{i+1})$ for all $1\leq i< n$. The
concatenation of these sequences without repeated overlaps is denoted
$\overline{u_1\cdots u_n}\in\mathcal{X}^{s+n-1}$.

The first representation we will consider is the circuit process of
MacQueen \cite{MacQueen1981df}. Let a \textit{circuit} be a periodic
function on $\mathcal{X}$, and consider a class of positively weighted
circuits $\mathscr{C}$ (for an example, see Figure~\ref{circuit process picture}).

\begin{defn}
A \textit{circuit process} of order $r$ is a Markov chain of the same
order, where the transition probability from $u\in\mathcal{X}^r$ to
any $v\in\mathcal{X}^r$ with $\Omega(u) = \mathrm{A}(v)$ is given by
\[
\frac{\sum_{\gamma\in\mathscr{C}} w_\gamma J_\gamma(\overline{uv}) }
{\sum_{\gamma\in\mathscr{C}} w_\gamma J_\gamma(u)},
\]
where $w_\gamma>0$ is the weight of circuit $\gamma$, and the function
$J_\gamma(\cdot)$ counts the number of times that the circuit traverses
a sequence in one period. In other words, in each step we move along
some circuit in $\mathscr{C}$ containing the current state with
probability proportional to its weight. The process only visits states
that appear in the circuits, for which transition probabilities are
well defined.
\end{defn}

An irreducible order-$r$ Markov chain with stationary law $P_{\pi}$ is
parametri\-zed by $P_\pi(u)$ for all $u\in\mathcal{X}^{r+1}$. One can
check that in a circuit process, this~is just $P_\pi(u) = \sum_{\gamma
\in\mathscr{C}} w_\gamma J_\gamma(u)$. MacQueen showed that any
order-$r$ Markov chain can be represented as a circuit process on a
finite set $\mathscr{C}$, which is not unique \cite{MacQueen1981df}.
This is true in particular when the chain is reversible.

%
\begin{figure}

\includegraphics{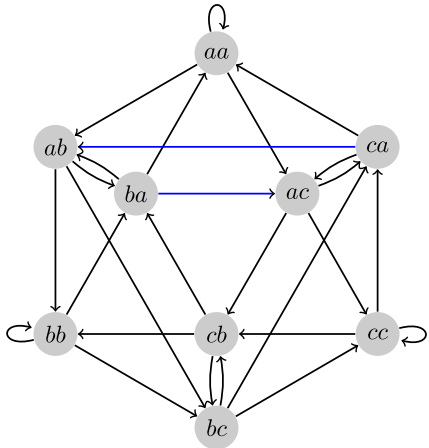}

\caption{A de Bruijn graph of order 2 on the state space $\mathcal{X}=\{
a,b,c\}$. In a reversible random walk, the two highlighted edges have
the same weight.}
\label{triangle}
\end{figure}

We introduce a second graphical representation that is canonical,
unlike the circuit process. Consider a \textit{de Bruijn graph} on the
vertices $\mathcal{X}^r$, which has a directed edge from $u$ to $v$ if
and only if $\Omega(u) = \mathrm{A}(v)$. That is, every path on the
graph is an admissible path. For an example, see Figure~\ref{triangle}.
Assign a weight $k_{uv}\geq0$ to each edge, and let $k_u$ be the
summed weights of edges departing from $u$. Furthermore, require
that
%
\begin{eqnarray}\label{eq1}
k_{uv} & =& k_{v^*u^*}  \qquad  \mbox{for every edge } uv, \\\label{eq2}
k_{u} & =& k_{u^*}  \qquad  \mbox{for all } u\in\mathcal{X}^r  \quad \mbox{and} \\\label{eq3}
\sum_{u\in\mathcal{X}^r} k_u &=& 1.
\end{eqnarray}

\begin{defn}
The \textit{reversible random walk of order $r$} is a random~walk on such
a graph, with transition probabilities
\[
p(v|u) = \frac{k_{uv}}{k_u}.
\]
\end{defn}

\begin{prop}
\label{representation}
An irreducible, reversible random walk of order $r$ represents a
reversible Markov chain of the same order. Every irreducible,
reversible order-$r$ Markov chain is equivalent to a unique reversible
random walk of order $r$.
\end{prop}

\begin{pf}
Let $\pi$ be the stationary distribution of the random walk. To prove
the first statement, we will first verify that $\pi(u)=k_u$ for all
$u\in\mathcal{X}^r$. Let $p(u|v)$ be the transition probability from
$v$ to $u$ in the random walk,~and recall that $\Omega(u) = \mathrm
{A}(v)$ iff $\Omega(v^*)=\mathrm{A}(u^*)$, then
\begin{eqnarray*}
\sum_{u\in\mathcal{X}^r} \pi(u)p(v|u) & =&
\sum_{\{u\in\mathcal{X}^{r}\dvtx \Omega(u)=\mathrm{A}(v)\}} k_u \frac
{k_{uv}}{k_u} \\
& =& \sum_{\{u\in\mathcal{X}^{r}\dvtx \Omega(v^*)=\mathrm{A}(u^*)\}}
k_{v^*u^*} = k_{v^*} = k_{v} = \pi(v).
\end{eqnarray*}
Then, the stationary law $P_\pi$ in the random walk of a path $u,v$ is
just
\[
P_\pi(u,v) = \pi(u)p(v|u) = k_u \frac{k_{uv}}{k_u} = k_{uv},
\]
which implies that $P_\pi(u,v) = k_{uv} = k_{v^*u^*} = P_\pi(v^*,u^*)$.
Therefore, the $\mathcal{X}$-va\-lued, order-$r$ Markov chain represented
by the random walk satisfies the reversibility condition in Proposition
\ref{ReversibilityCondition}. Proving the second statement is now
straightforward. Let $V_n,n\in\mathbb{N}$, be the first-order
representation of an irreducible, order-$r$ Markov chain, with
transition probabilities $p'(v|u)$. By the Perron--Frobenius theorem,
$V$ has a unique stationary distribution $\pi'$. Assign edge weights to
the de Bruijn graph on $\mathcal{X}^r$, setting $k_{uv} = \pi
'(u)p'(v|u)$. Since the order-$r$ Markov chain is reversible, it
follows directly from Proposition \ref{ReversibilityCondition} that the
edge weights satisfy conditions \eqref{eq1}--\eqref{eq3}.
 \end{pf}

\section{From a reinforced random walk to the conjugate prior}\label{sec3}

An edge-re\-inforced random walk (ERRW) is a random walk on an finite,
undirected graph, where every edge-weight is increased by 1 each time
it is crossed. Since Diaconis and Coppersmith defined this process
\cite{Diaconis1988}, we have learned that it is partially
exchangeable and, by de Finetti's theorem for Markov chains, a~mixture
of Markov chains \cite{Diaconis1980hb}. The mixing measure, which
lives on the space of reversible Markov chains, was more recently
characterized in the literature \cite{Keane2000df}. Diaconis and
Rolles showed that this distribution is a conjugate prior for the
reversible Markov chain, much as the Beta distribution, arising from a
Poly\'{a} urn scheme, is a conjugate prior for sequences of i.i.d. binary
random variables \cite{Diaconis2006gd}.

Here, we construct a conjugate prior for higher-order reversible Markov
chains via a reinforced random walk in $\mathcal{X}^r$, making use of
de Finetti's theorem for Markov chains. This process is markedly
different from an ERRW in $\mathcal{X}^r$ due to the structure of a
reversible Markov chain with memory, although it is designed to be
partially exchangeable.

Let $\alpha$ be any sequence on $\mathcal{X}$ and $v$ a sequence
shorter than $\alpha$. Define the function $J'_\alpha(v)$, which counts
the number of times that $v$ appears in $\alpha$, and $J''_\alpha(v)$,
which counts the number of times that $v$ appears in $\alpha$ followed
by at least one state. Fix $w$, a stationary measure for an
irreducible, reversible, order-$r$ Markov chain. Also fix $v_0\in
\mathcal{X}^r$. Let $\beta$ be a palindromic sequence that starts with
$v_0$ and ends with $v_0^*$. Choose a positive constant $c$, such that
for all $u\in\mathcal{X}^{r+1}$, $w(u)-c J'_\beta(u) > 0$. Now, given a
sequence $\eta$, starting with $v_0$, and any sequence $v$, define the functions
%
\begin{eqnarray}
\label{w'}
w'(\eta,v) & =& w(v) + c \bigl( J'_\eta(v) + J'_{\eta^*}(v) - J'_\beta(v)\bigr)
  \quad \mbox{and} \\
\label{w''}
w''(\eta,v) & =& w(v) + c \bigl( J''_\eta(v) + J''_{\eta^*}(v) - J''_\beta(v)\bigr).
\end{eqnarray}
When $\eta$ represents the path of a stochastic process in $\mathcal
{X}^r$ up to time $n$ (formally, $\eta=\overline{v_0\cdots v_n}$), we
will use the notation $w'_n(v) \equiv w'(\eta,v)$ and $w''_n(v) \equiv
w''(\eta,v)$.

\begin{defn}
\label{reinforcement}
The \textit{reinforced random walk of order $r$} is a stochastic process
$Y_n, n\in\mathbb{N}$, on $\mathcal{X}^r$ with distribution $Q_{w,v_0}$.
The initial state is $v_0$ with probability 1. For any admissible path
$v_0,\ldots ,v_n$, the conditional transition probability
\[
Q_{w,v_0}(Y_{n+1}=u | Y_0=v_0,\ldots ,Y_n=v_n) = \frac{w'_n(\overline
{v_n u})}{w''_n(v_n)}
\]
whenever $v_n,u$ is admissible and zero otherwise.
\end{defn}

\begin{rmk}
The law $Q_{w,v_0}$ also depends on $\beta$ and $c$. These parameters
are constant in the following discussion, so they are omitted from the
notation for conciseness. When $r=1$, this process is equivalent to an
ERRW. In this case, the palindrome is unnecessary because the terms
involving $\beta$ in $w'$ and $w''$ can be modeled with a different
$w$. For $r\geq2$, this is not the case, and $\beta$ is essential for
partial exchangeability (see Proposition \ref{PE}).
\end{rmk}

\begin{rmk}
\label{circuit_interpretation}
This process admits an interpretation as a reinforcement scheme of the
circuit process. Consider a circuit process of order $r$ with
stationary probability $w(u)=\sum_{\gamma\in\mathscr{C}}w_\gamma
J_\gamma(u)$ for all $u\in\mathcal{X}^{r+1}$. In addition, consider
three weighted sequences: the palindrome $\beta$, a sequence $\eta$
that represents the path of the reinforced process from the initial
state $v_0$ up to the current state, and the reversed path $\eta^*$.
These are depicted in Figure~\ref{reinforcement_picture} along with
their weights $-c$, $c$  and $c$, respectively. As in the circuit
process, we move along any circuit or sequence that contains the
current state with probability proportional to its weight. The
reinforcement is accomplished by elongating the paths $\eta$ and $\eta^*$.
\end{rmk}

\begin{rmk}
\label{errw_interpretation}
The process is also a reinforcement scheme of a modified reversible
random walk of order $r$. Consider a weighted de Bruijn graph, where
for every admissible $u,v$, $k_{uv} = w(\overline{uv})$. Then, for
every $\overline{uv}$ in the palindrome $\beta$, subtract $c$ from
$k_{uv}$. The reinforcement scheme will consist of a random walk on the
resulting graph, where after every transition $v_{i}\to v_{i+1}$ we
increase both $k_{v_iv_{i+1}}$ and $k_{v_{i+1}^*v_i^*}$ by $c$.
Accordingly, if $\overline{v_iv_{i+1}}$ is a~palindrome,\vspace*{2pt} the weight
$k_{v_iv_{i+1}}$ is increased by $2c$.\vspace{-6pt}
\end{rmk}

\begin{prop}
\label{PE}
The reinforced random walk of order $r$ is partially~ex\-changeable in
the sense of Diaconis and Freedman \cite{Diaconis1980hb}.
\end{prop}

\begin{pf}
We must show that the probability $Q_{w,v_0}(v_0,\ldots ,v_n)$ of any
admissible path $v_0,\ldots ,v_n$ is a function of the initial state
$v_0$ and the transition counts between every pair of states. For any
pair $u,v$ in $\mathcal{X}^r$ with $\mathrm{A}(v)=\Omega(u)$, let
$C(u,v)$ be the total number of transitions $u\to v$, and $v^*\to u^*$.
We will show the stronger statement that $v_0$ and $C$ are sufficient
statistics for the reinforced random walk.
%
%
\begin{figure}

\includegraphics{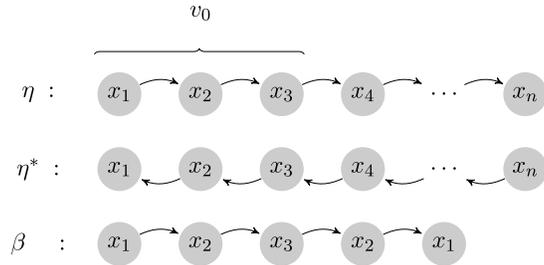}

\caption{Auxiliary sequences in the order-$r$ reinforced random walk.}
\label{reinforcement_picture}
\end{figure}

Let us first establish some properties that are conserved in the
process. For every $u\in\mathcal{X}^{r+1}$, the initial weights
$w'_0(u)$ and $w'_0(u^*)$ are equal. This is direct from the definition
in equation~\eqref{w'} because: $w$ defines a reversible Markov chain of
order $r$; the functions $J_{v_0}'$ and $J_{v_0^*}'$ are zero for both
$u$ and $u^*$; and $\beta$ is a palindrome, so if it contains $u$, it
also contains $u^*$, and $J_\beta'(u)=J_\beta'(u^*)$. This property is
maintained after every transition $v_n\to v_{n+1}$, because the weights
may both be increased by $c$ if $\overline{v_nv_{n+1}}$ is $u$ or
$u^*$, or both remain constant otherwise.

For every $v\neq v_0$ in $\mathcal{X}^r$, the initial weights
$w''_0(v)=w''_0(v^*)$. This is direct from equation~\eqref{w''} because: $w$
is reversible, both $J_{v_0}''$ and $J_{v_0^*}''$ are zero for $v$ and
$v^*$, and the sequence $\beta$ is a palindrome, so for every
transition starting at $v$ there will be another starting from $v^*$.
The last fact is not necessarily true for $v_0$, because unless $v_0$
itself is a palindrome, $\beta$ will contain a transition starting from
it, but no transition starting from $v_0^*$. So, in the beginning,
$w''_0(v_0) = w''_0(v_0^*)-c$. When a transition occurs from $v_0$ to
$v_1$, the weights $w''_1(v_0)$ and $w''_1(v_0^*)$ become equal, while
$w''_1(v_1) = w''_1(v_1^*)-c$, provided $v_1$ is not a palindrome.
Hence, this singularity is preserved by the last state visited by the process.

The probability $Q_{w,v_0}(v_0,\ldots,  v_n)$ is a ratio of two products.
In the numerator, we find a factor of the form $w'_t(\overline{uv})$
for every admissible transition $u\to v$, while in the denominator, we
find a corresponding weight $w''_t(u)$. It is easy to check that the
numerator is only a function of $C$. Every transition $u\to v$ or
$v^*\to u^*$ adds a new factor of $w'_t(\overline{uv})$, which is
always greater than the previous one by $c$. If $\overline{uv}$ is a
palindrome, then every new factor of $w'_t(\overline{uv})$ is increased
by $2c$. So, the numerator can be computed from the initial weights and $C$.

We have left to show that the denominator is only dependent on $v_0$
and $C$. Note that the transition counts from $v$ or $v^*$ are a
function of $C$ and $v_0$, because every event $v\to u$ is a transition
from $v$, while every event $u^*\to v^*$ is followed by a transition
from $v^*$, unless this is the final state, which is determined by
$v_0$. After every transition from $v$ or $v^*$, we add a factor of
$w_t''(v)$ or $w''_t(v^*)$ to the denominator. At any time $t$, these
weights differ by $c$ (if $v$ is not a palindrome), but the factor
added is always the smaller of the two. Between two transitions, each
of these weights is reinforced by $c$, so consecutive factors differ by
that amount. If $v$ is a palindrome, there is no distinction between
$w''_t(v)$ and $w''_t(v^*)$, and consecutive factors differ
by~$2c$.
\end{pf}

\begin{lem}
\label{positive-recurrence}
Suppose that in the reinforced random walk, we visit $v$ and $v^*$ in
$\mathcal{X}^r$ infinitely often a.s., and let $\tau_n$ be the $n$th
time we visit either state. The process $Y_{\tau_n}$ is a mixture of
Markov chains. Furthermore, if $D_n$ is the ratio of the number of
visits to $v^*$ and $v$ by $\tau_n$, $D_n$ converges a.s. to a finite
limit $D_\infty$.
\end{lem}

\begin{pf}
We claim that if $Y_n$ is partially exchangeable, so is $Y_{\tau_n}$.
It is sufficient to show that the probability of a sequence $Y_{\tau
_n}$ is invariant upon block transpositions, which generate the group
of permutations that preserve transition counts (\cite{Diaconis1980hb},
Proposition~27). The probability of a path $v_{\tau_1},
\ldots
,v_{\tau_n}$ in $Y_{\tau_n}$ is the sum of the probabilities of all
paths $v_0,v_1,\ldots ,v_{\tau_n}$ in $Y_n$ that map to it. Denote this
set of paths $\Theta$. After a transposition of $v$-blocks or
$v^*$-blocks, the probability of the path in $Y_{\tau_n}$ is equal to
the sum of the probabilities of a different set of paths $\Theta'$ in
$Y_n$. However, it is easy to see that this transpostion of $v$-blocks
or $v^*$-blocks defines a bijection from $\Theta$ to $\Theta'$, and the
probability of each path and its transposition is the same, because
$Y_n$ is partially exchangeable. Therefore, $Y_{\tau_n}$ is partially
exchangeable. Furthermore, we assume that $v$ and $v^*$ are recurrent,
so by de Finetti's theorem for Markov chains $Y_{\tau_n}$ is a mixture
of Markov chains with a unique measure $\mu$ on the space of 2 by 2
transition matrices \cite{Diaconis1980hb}. Note that both states are
recurrent with probability 1, so the subset of transition matrices
where one of the states is transient has $\mu$-measure zero. This
implies that $\mu$-a.s.  the transition matrix is irreducible, and
since the state space is finite, both states are positive-recurrent.
Therefore, $D_n$ converges a.s.  to a finite limit.\looseness=-1
\end{pf}

\begin{prop}
The reinforced random walk of order $r$ traverses~eve\-ry edge $v \to u$
with $w(\overline{uv})>0$ infinitely often, almost surely.
\label{recurrence}
\end{prop}

\begin{pf}
As $\mathcal{X}$ is finite, we must visit at least one state in
$\mathcal{X}^r$ infinitely often, so without loss of generality, let
this state be $v$. Let $\tau_n$ be the $n$th time we visit $v$, and
$\mathcal{F}_n$ be $\sigma(Y_1,\ldots ,Y_{\tau_n})$. For $u$ with $v,u$
admissible and $w(\overline{vu})>0$, let $A_n$ be the event that
$Y_{\tau_n +1}= u$. Also, let $p_n = Q_{w,v_0}(A_n|\mathcal{F}_n)$. By
L\'{e}vy's extension of the Borel--Cantelli lemma (Lemma
\ref{levy}),
\[
\lim_{n\to\infty}\frac{ \sum_{m=1}^n 1_{A_m} }{ \sum_{m=1}^n p_m }= 1
  \qquad \mbox{on }   \Biggl\{ \sum_{m=1}^\infty p_m = \infty \Biggr\}.
\]
Therefore, to show that the transition $v\to u$ is observed infinitely
often with probability 1, it is sufficient to show that $\sum_m p_m =
\infty$ a.s. The conditional probability $p_m$ is just $w'_{\tau
_m}(\overline{vu})/w_{\tau_m}''(v)$. Let $B_{m,k}$ be the event that we
observe $v^*$ fewer than $km$ times between $\tau_1$ and $\tau_m$. On
$B_{m,k}$, we can lower-bound $p_m$ using the minimum possible value of
$w'_{\tau_m}(\overline{vu})$, which is its initial value, and the
maximum possible value of $w''_{\tau_m}(v)$, which is $(k+1)mc$.
Thus,
\begin{eqnarray*}
p_m &=& Q_{w,v_0}(A_m\cap B_{m,k}|\mathcal{F}_n)+ Q_{w,v_0}(A_m\cap
B_{m,k}^C|\mathcal{F}_n) \\
&\geq&\mathbf{1}_{B_{m,k}}\frac{w'_{\tau_1}(\overline{vu})}{(k+1)mc}.
\end{eqnarray*}
Now, consider the event $\{D_\infty< N\}$. On this set, for any $k>N$,
we will be in $B_{m,k}$ for all but finitely many $m$, which implies
$\sum_m p_m = \infty$, by the previous inequality. But, by Lemma \ref
{positive-recurrence} we have $Q_{w,v_0}\{D_\infty<\infty\}=1$, so
noting $\{D_\infty< \infty\}=\bigcup_{N\in\mathbb{N}}\{D_\infty<N\}$ we
conclude that $\sum_m p_m =\infty$ $Q_{w,v_0}$-a.s., and $A_m$ happens
infinitely often. Since $w$ defines an irreducible Markov chain, the
proposition follows by induction.
\end{pf}

Propositions \ref{recurrence} and \ref{PE} are sufficient to show by de
Finetti's theorem for Markov chains \cite{Diaconis1980hb} that the
reinforced random walk of order $r$ is a mixture of Markov chains on
$\mathcal{X}^r$, or
%
\begin{equation}
\label{mixture}
Q_{w,v_0}(v_0,\ldots ,v_n) = \int_\mathcal{T} P_{v_0}^T(v_0,\ldots
,v_n)\,d\phi_{w,v_0}(T),
\end{equation}
where $P_{v_0}^T$ is the distribution of a Markov chain started at
$v_0$ and parametri\-zed by the matrix $T$, $\mathcal{T}$ is the space of
$\mathcal{X}^r\times\mathcal{X}^r$ stochastic matrices  and $\phi
_{w,v_0}$ is a unique measure on the Borel subsets of this space. Let
$\mathcal{T'}\subseteq\mathcal{T}$ be the set of matrices that
represent irreducible, reversible Markov chains of order $r$.

\begin{prop}
The reinforced random walk of order $r$ is a mixture of reversible
Markov chains of the same order, or $\phi_{w,v_0}(\mathcal{T}')=1$.
\end{prop}

\begin{pf}
This is a special case of Proposition \ref{Concentration}.
\end{pf}

\section{Variable-order, reversible Markov chains}\label{sec4}

The number of parameters of a Markov chain grows as $|\mathcal{X}|^r$
with the order, $r$, which renders higher-order models impractical in
many statistical applications. In this section, we investigate a family
of models with finite memory length which do not suffer from this curse
of dimensionality.

\begin{defn}
\label{mixed-order definition}
A \textit{variable-order Markov chain} is a Markov chain of order $r$
with the constraint that for every \textit{history} $h$ in the set
$\mathscr{H}\subseteq\{v\in\mathcal{X}^q\dvtx q< r \}$, if two states
$u,u'\in\mathcal{X}^r$ both end in $h$, the transition probabilities
$p(v|u)$ and $p(v|u')$ are equal for every $v\in\mathcal{X}^r$.
\end{defn}

In essence, this is a discrete process which upon reaching a sequence
$h\in\mathscr{H}$ loses memory of what preceded it. When $\mathscr{H}$
is empty, we recover a~general Markov chain of order $r$.
Variable-order Markov chains have proven useful in applications where
there is long memory only in certain directions. The literature on the
subject can be traced to Rissanen \cite{Rissanen1983} and Weinberger
\cite{Weinberger1995}, who developed tree-based algorithms for
estimating the set of histories efficiently in the context of
compression. B\"{u}hlmann and Wyner proved several consistency results
on these algorithms \cite{Buhlmann1999}, and the former later
addressed the problem of model selection \cite{Buhlmann2000}. For an
evaluation of different algorithms in applications, see \cite{Begleiter2004}.

It is worth noting that MacQueen mentioned variable-order Markov chains
in an unpublished abstract. However, there is a marked difference
between his definition and B\"{u}hlmann and Wyner's, which relates to
the closure properties of $\mathscr{H}$. MacQueen requires that if $h$
is in $\mathscr{H}$, then so are all the sequences that begin with~$h$.
Intuitively, this means that the process cannot recover memory once it
is lost. B\"{u}hlmann and Wyner do not impose this constraint. However,
this is guaranteed when the process is reversible.

\begin{prop}
\label{regaining_memory}
Let $X_n, n\in\mathbb{N}$, be an irreducible, reversible, variable-order
Markov chain with histories $\mathscr{H}$. If $h\in\mathscr{H}$, then
$h^*$ is also a history; additionally, any sequence that has $h$ as a
prefix is also in $\mathscr{H}$.
\end{prop}

\begin{pf}
Let $P_\pi$ be the stationary law of the chain. If $h\in\mathscr{H}$,
then for any pair $a,b\in\mathcal{X}^q$, where $q$ and the length of
$h$ sum to $r$, $P_\pi( X_1,\ldots , X_{r+q}= ahb | X_1,\ldots , X_r = ah)$
is independent of $a$, or
\[
\frac{P_\pi(ahb)}{P_\pi(ah)} = C   \qquad  \forall a\in\mathcal{X}^q.
\]
This implies
\[
\frac{P_\pi(hb)}{P_\pi(h)} = \frac{\sum_{a\in\mathcal{X}^q} P_\pi
(ahb)}{\sum_{a\in\mathcal{X}^q} P_\pi(ah)}
= \frac{\sum_{a\in\mathcal{X}^q} P_\pi(ah) C}{\sum_{a\in\mathcal{X}^q}
P_\pi(ah)} = \frac{P_\pi(ahb)}{P_\pi(ah)}.
\]
Using the fact that $P_\pi$ is invariant upon time reversal and
rearranging factors, we obtain
\begin{eqnarray*}
\frac{P_\pi(b^*h^*a^*)}{P_\pi(b^*h^*)} = \frac{P_\pi(h^*a^*)}{P_\pi(h^*)}.
\end{eqnarray*}
The left-hand side is equal to $P_\pi( X_1,\ldots , X_{r+q}=b^*h^*a^*|
X_1,\ldots ,X_r = b^*h^*)$, which by the previous identity is independent
of $b^*$. As this is true for any $a\in\mathcal{X}^q$, $h^*$ must be a
history in $\mathscr{H}$. To prove the second part of the statement,
suppose $h$ is a prefix of $g$. Since $h^*$ is in $\mathscr{H}$, and
$g^*$ ends in $h^*$, then by definition $g^*\in\mathscr{H}$. Using the
first result, we conclude that $g\in\mathscr{H}$.
\end{pf}

We will define a reinforcement scheme, which like the one in the
previous section is recurrent, partially exchangeable and, by de
Finetti's theorem, a~mixture of Markov chains. But, in this case, the
mixing measure is restricted to the variable-order, reversible Markov
chains with a fixed set of histories $\mathscr{H}$. As before, we begin
with a stationary, reversible function $w$, an initial state $v_0\in
\mathcal{X}^r$, and a palindromic sequence $\beta$ that starts with
$v_0$. Let the function $f\dvtx \mathcal{X}^r\mapsto\mathscr{H}$ map any
sequence to its shortest ending in $\mathscr{H}$.

\begin{defn}
\label{mixed-order reinforcement}
The \textit{variable-order}, \textit{reinforced random walk} is a stochastic
process $Z_n, n\in\mathbb{N}$, on $\mathcal{X}^r$ with measure
$H_{w,v_0}$. The initial state is $v_0$ with probability 1. For any
admissible path $v_0,\ldots ,v_n$, the conditional transition probability
\[
H_{w,v_0}(Z_{n+1}=u | Z_0=v_0,\ldots ,Z_n=v_n) = \frac{w'_n(\overline
{f(v_n) u})}{w''_n(f(v_n))}
\]
whenever $v_n,u$ is admissible and zero otherwise.
\end{defn}

\begin{rmk}
This process is a reinforced circuit process, just like the one defined
in Remark \ref{circuit_interpretation}, with the difference that in
computing the transition probabilities, instead of taking the current
state to be the sequence $v_n\in\mathcal{X}^r$, we let it be the
shortest ending of $v_n$ in $\mathscr{H}$, or $f(v_n)$.
\end{rmk}

\begin{prop}
The variable-order, reinforced random walk is partially exchangeable in
the sense of Diaconis and Freedman.
\label{mixed-PE}
\end{prop}

This proof is deferred to the \hyperref[sec1]{Appendix}. One can show that this process
is recurrent following the same argument of Proposition \ref
{recurrence}. In the proof of Proposition~\ref{recurrence}, we use a
shortest history $h$ in place of $v$, and Lemma \ref
{positive-recurrence} still holds for $h$ and $h^*$. Recurrence and
partial exchangeability imply
%
\begin{equation}
\label{mixed-order mixture}
H_{w,v_0}(v_0,\ldots ,v_n) = \int_\mathcal{T} P_{v_0}^T(v_0,\ldots
,v_n)\,d\psi_{w,v_0}(T)
\end{equation}
for a unique measure $\psi_{w,v_0}$ characterized by the function $w$,
and the initial state, in addition to the parameters $\beta$, $c$  and
$\mathscr{H}$, which we keep fixed. In the \hyperref[sec1]{Appendix}, we show that $\psi
_{w,v_0}$ is restricted to the reversible, variable-order Markov chains
with histories $\mathscr{H}$.

\begin{prop}
\label{Concentration}
Let $\mathcal{T}''\subseteq\mathcal{T}$ be the set of transition
matrices representing an irreducible, reversible, variable-order Markov
chain where every $h\in\mathscr{H}$ is a history. Then, $\psi
_{w,v_0}(\mathcal{T}'')=1$.
\end{prop}

\section{Bayesian analysis}\label{sec5}

In Section~\ref{sec3}, we defined a family of measures in the space of
order-$r$, reversible Markov chains, and in Section~\ref{sec4} we extended it to
variable-order, reversible Markov chains. In the following, we will
show that these distributions are conjugate priors for a Markov chain
of order~$r$. We discuss properties of the prior relevant to Bayesian
analysis, such as a~natural sampling algorithm and closed-form
expressions for some important moments.

\begin{defn}
Consider a variable-order, reinforced random walk $Z_n$,
$n\in\mathbb{N}$,
with distribution $H_{w,v_0}$ and take any admissible path $e= v_0,\ldots
,v_n$. We define $Z^{(e)}_n, n\in\mathbb{N}$, to be the process with
law
\begin{eqnarray*}
 &&H_{w,v_0,e}(v_n,u_1,\ldots ,u_m) \\
&& \qquad  =  H_{w,v_0}(Z_{n+1}=u_1,\ldots ,Z_{n+m+1}=u_m | Z_1=v_1,\ldots , Z_m=v_m).
\end{eqnarray*}
\end{defn}

In words, $Z^{(e)}$ is the continuation of a variable-order reinforced
random walk after traversing some fixed path $e$. We can rewrite the law
%
\begin{equation}
\label{post}
H_{w,v_0,e}(v_n,u_1,u_2,\ldots ,u_m) = \frac{H_{w,v_0}(v_1,\ldots
,v_n,u_1,\ldots ,u_m)}{H_{w,v_0}(v_1,\ldots ,v_n)},
\end{equation}
which makes it evident that $Z^{(e)}$ is partially exchangeable,
because for a~fi\-xed~$e$, the numerator only depends on the transition
counts in $v_n,u_1,\ldots ,u_m$, while the denominator is constant. It is
also not hard to see that the process visits every state infinitely
often with probability 1. Therefore, by de Finetti's theorem for Markov
chains, it is a mixture of Markov chains with a mixing measure that
will be denoted $\psi_{w,v_0,e}$.

\begin{prop}
Suppose we model a process $W_n, n\in\mathbb{N}$, as a reversible,
variable-order Markov chain with histories $\mathscr{H}\subseteq\{v\in
\mathcal{X}^q\dvtx q< r \}$, and we assign a prior $\psi_{w,v_0}$ to the
transition probabilities, $T$. Given an observed path, $e=v_0,\ldots ,v_n$,
the posterior probability of $T$ is $\psi_{w,v_0,e}$. In consequence,
the family of measures
\[
\mathcal{D}=\{\psi_{w,v_0,e}\dvtx e \mbox{ an admissible path starting in
}v_0\}
\]
is closed under sampling.
\end{prop}

\begin{pf}
Consider the event $W_n=v_n,W_{n+1}=u_1,\ldots ,W_{n+1+m}=u_m$. By Bayes
rule, the posterior probability of this event given the observation is
the prior probability of $W_{1}=v_1,\ldots ,W_n=v_n,W_{n+1}=u_1,\ldots
,W_{n+1+m}=u_m$ divided by the prior probability of $W_{1}=v_1,\ldots
,W_n=v_n$. By equation~\eqref{post}, this posterior is equal to
$H_{w,v_0,e}$. Let $\rho(T)$ be the posterior distribution of~$T$ given
the observation, then for any $u_1,\ldots ,u_m$ and any $m>0$,
\[
H_{w,v_0,e}(v_n,u_1,\ldots ,u_m) = \int_\mathcal{T}
P_{v_n}^T(v_n,u_1,\ldots ,u_m)\,d\rho(T).
\]
By de Finetti's theorem for Markov chains, the mixing measure $\psi
_{w,v_0,e}$ is unique; therefore, we must have $\rho=\psi_{w,v_0,e}$.
\end{pf}

In the next proposition, we show that the variable-order, reinforced
random walk may be used to simulate from the conjugate prior $\psi
_{w,v_0}$ (or~\mbox{using} a similar argument, a posterior of the form $\psi
_{w,v_0,e}$). Let $\{V^{(i)} = v_1^{(i)},v_2^{(i)},\ldots ,\break v_n^{(i)}\}
_{i\in\{1,\ldots ,k\}}$ be independent samples of the reinforced random
walk with initial parameters $w$ and $v_0$. For any sequence $u\in
\mathcal{X}^{r+1}$, consider the random variable $n^{-1}w'_n(u)$, the
weight defined in equation~\eqref{w'} for a sample path with distribution
$H_{w,v_0}$, normalized by the path's length. Define the empirical
estimate, $n^{-1}w'_{n,k}(u)$, to be the mean of this random variable
evaluated at the paths $\{V^{(i)}\}_{i\in\{1,\ldots ,k\}}$. Also, let
$P^T_\pi$ be the stationary law of an order-$r$ Markov chain with
transition probabilities $T$. We have seen that $\{P_\pi^T(u)\dvtx u\in
\mathcal{X}^{r+1}\}$ has a one-to-one correspondence with $T$.

\begin{prop} For any bounded, real-valued function $g(P_\pi^T(\cdot))$,
\label{posterior_sampling}
%
\begin{equation}
\label{simulation}
\lim_{n\to\infty}\lim_{k\to\infty} g\bigl(\{n^{-1}w'_{n,k}(u)\dvtx u\in\mathcal
{X}^{r+1}\}\bigr) \stackrel{a.s.}{=} \int_\mathcal{T}
g(P_\pi^T)\,d\psi_{w,v_0}(T).
\end{equation}
\end{prop}

\begin{pf}
The empirical estimate $g(\{n^{-1}w'_{n,k}(u)\dvtx u\in\mathcal{X}^{r+1}\})$
is the average of i.i.d. observations, so by the strong law of large
numbers, w.p.1,
\[
\lim_{k\to\infty} g\bigl(\{n^{-1}w'_{n,k}(u)\dvtx u\in\mathcal{X}^{r+1}\}\bigr) =
H_{w,v_0} \bigl[ g\bigl(\{n^{-1}w'_{n}(u)\dvtx u\in\mathcal{X}^{r+1}\}\bigr)  \bigr],
\]
where the right-hand side is the expectation in a reinforced random
walk with parameters $w,v_0$. In the proof of Proposition \ref
{Concentration}, we showed that $w'_n(u)$ converges $H_{w,v_0}$-a.s.
Taking the limit as $n\to\infty$, by dominated
convergence,
\begin{eqnarray*}
&& \lim_{n\to\infty}\lim_{k\to\infty} g\bigl(\{n^{-1}w'_{n,k}(u)\dvtx u\in
\mathcal{X}^{r+1}\}\bigr)  \\
&& \qquad  = H_{w,v_0} \Bigl[ \lim_{n\to\infty} g\bigl(\{n^{-1}w'_{n}(u)\dvtx u\in\mathcal
{X}^{r+1}\}\bigr)  \Bigr].
\end{eqnarray*}
Conditional on a variable $T$ measurable on its tail $\sigma$-field
with distribution $\psi_{w,v_0}$, the reinforced random walk is a
Markov chain with law $P_{v_0}^T$. We know $w'_n(u)$ converges
$P_{v_0}^T$-a.s.  to $P_\pi^T(u)$, so equation~\eqref{simulation} follows.
\end{pf}

Several moments of $H_{w,v_0}$ have closed-form expressions. In
particular, the mean likelihood $P_{v_0}^T$ of any path beginning in
$v_0$ is just the probability of the path in the reinforced random walk
by equation~\eqref{mixed-order mixture}. From the proof of Proposition \ref
{mixed-PE}, one can deduce a closed-form expression for the law of the
variable-order reinforced random walk as a function of the transition
counts in a path (see Supplement \cite{Bacallado2010}). From a
realization of the transition counts as a path, one can also compute
the law $H_{w,v_0}$ by modeling a random walk with reinforcement.

The expectation of cycle probabilities with a prior $\psi_{w,v_0}$ on
$T$ may also be computed exactly.

\begin{prop}
For any cyclic path $v,v_1,\ldots ,v_n,v$, not necessarily including
$v_0$, the expectation of $P_{v}^T(v,v_1,\ldots ,v_n,v)$ with prior $\psi
_{w,v_0}$ on $T$ has a closed-form expression, provided $w'_0(u)$ is
greater than $3c$ for all $u\in\mathcal{X}^{r+1}$.
\end{prop}

\begin{pf}
Find the shortest cycle $v,\ldots ,v_0,\ldots ,v$ with positive weight
$w$. Then, for any transition matrix $T$ in the support of $\psi
_{w,v_0}$, we have
%
\begin{equation}
\label{ratio}
P_{v}^T(v,v_1,\ldots ,v_n,v) = \frac{P_{v_0}^T(v,v_1,\ldots ,v_n,v,\ldots
,v_0,\ldots ,v)}{P_{v_0}^T(v,\ldots ,v_0,\ldots ,v)}.
\end{equation}
Taking the expectation with a measure $\psi_{w,v_0}$ on $T$, we obtain
\begin{eqnarray*}
 && \int_{\mathcal{T}} P_{v}^T(v,v_1,\ldots ,v_n,v)\,d\psi
_{w,v_0}(T)   \\
&& \qquad  = \int_{\mathcal{T}} \frac{P_{v_0}^T(v,v_1,\ldots ,v_n,v,\ldots
,v_0,\ldots ,v)}{P_{v_0}^T(v,\ldots ,v_0,\ldots ,v)}\,d\psi_{w,v_0}(T).
\end{eqnarray*}
By Bayes theorem, the product of the likelihood $P_{v_0}^T(v,v_1,\ldots
,v_n,v,\ldots , v_0,\break\ldots ,v)$ and the prior $d\psi_{w,v_0}(T)$ is equal
to the marginal prior probability of the path $v,v_1,\ldots ,v_n,v,\ldots
,v_0,\ldots ,v$ times the posterior of $T$:
\begin{eqnarray*}
 && \int_{\mathcal{T}} P_{v}^T(v,v_1,\ldots ,v)\,d\psi_{w,v_0}(T)
  \\
&& \qquad = H_{w,v_0}(v,v_1,\ldots ,v,\ldots ,v_0,\ldots ,v) \int_{\mathcal{T}}
\frac{1}{P_{v_0}^T(v,\ldots ,v_0,\ldots ,v)}\,d\psi_{w_p,v_0}(T),
\end{eqnarray*}
where $w_p$ are the weights parametrizing the posterior of $T$ given
the path $v,v_1,\ldots ,v,\ldots ,v_0,\ldots ,v$. To solve the integral on
the right-hand side, let us rewrite it using Bayes theorem and equation~\eqref
{mixed-order mixture},
\begin{eqnarray*}
 &&  H^{-1}_{w_{pp},v_0}(v,\ldots ,v_0,\ldots ,v,\ldots ,v_0,\ldots
,v)   \\
&& \qquad {} \times\int_\mathcal{T} \frac{P_{v_0}^T(v,\ldots ,v_0,\ldots ,v,\ldots
,v_0,\ldots ,v)}{P_{v_0}^T(v,\ldots ,v_0,\ldots ,v)}\,d\psi_{w_{pp},v_0},
\end{eqnarray*}
where $w_{pp}$ are the weights $w_p$ reduced by the cycle $v,\ldots
,v_0,\ldots ,v,\ldots ,v_0,\ldots ,v$. These weights are positive because
of the assumption $w'_0(u)>3c$ for all $u$, which could certainly be
relaxed in some cases. Applying equations~\eqref{mixed-order mixture} and \eqref
{ratio} once more, the last expression becomes
\[
H^{-1}_{w_{pp},v_0}(v,\ldots ,v_0,\ldots ,v,\ldots ,v_0,\ldots
,v)H_{w_{pp},v_0}(v,\ldots ,v_0,\ldots ,v),
\]
which completes our derivation.
\end{pf}

The ability to compute these expectations exactly makes it possible to
use Bayes factors for model comparison \cite{Kass1995}. Given some
data $\mathbf{X}$ and two probabilistic models, where each model $i$
has a prior measure $P^{(i)}$ and parameters $\theta_i$, a Bayes factor
quantifies the relative odds between them. It is formally defined
as,\vspace{-5pt}
%
\begin{equation}
\label{bayes_factor}
\frac{P^{(1)}(\mathbf{X})}{P^{(2)}(\mathbf{X})}
= \frac{\int P^{(1)}(\mathbf{X}|\theta_1)\,dP^{(1)}(\theta_1)}{\int
P^{(2)}(\mathbf{X}|\theta_2)\,dP^{(2)}(\theta_2)},
\end{equation}
the ratio between the marginal probabilities of the data under each
model. Each marginal probability is sometimes referred to as the \textit
{evidence} for the corresponding model. Diaconis and Rolles apply Bayes
factors to compare a number of models on different data sets. They
consider reversible Markov chains, general Markov chains, and i.i.d.
models \cite{Diaconis2006gd}, assigning conjugate priors which
facilitate computing the marginal probabilities in
equation~\eqref{bayes_factor}.

The conjugate priors introduced here facilitate similar comparisons,
where the family of models under consideration is expanded to include
reversible Markov chains that differ in their length of memory. For
some data $\mathbf{X}$, one can define two variable-order reversible
Markov models, with different histories, $\mathscr{H}^{(1)}$ and
$\mathscr{H}^{(2)}$. In each case, we assign a conjugate prior, $\psi
_{w,v_0}^{(1)}$ and $\psi_{w,v_0}^{(2)}$, respectively, to the
transition probability matrix. To make the prior uninformative in some
sense we could set $w$ to be uniform for all $u\in\mathcal{X}^{r+1}$
and let $\beta$ be the shortest palindrome starting with $v_0$, for
example. The constant $c$ is set to 1. The Bayes factor is then\vspace{-2pt}
\[
\frac{P^{(1)}(\mathbf{X})}{P^{(2)}(\mathbf{X})}
=\frac{\int_{\mathcal{T}} P^T_{v_0}(\mathbf{X})\,d\psi
^{(1)}_{w,v_0}(T)}{\int_{\mathcal{T}} P_{v_0}^T(\mathbf{X})\,d\psi
^{(2)}_{w,v_0}(T)}.
\]

%
\begin{figure}[b]

\includegraphics{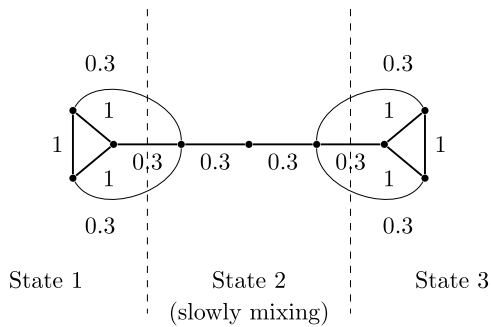}
\vspace{-3pt}
\caption{A lumped reversible Markov chain.}
\label{example}
\vspace{-5pt}
\end{figure}

We have seen that the expectations on the right-hand side can be
computed exactly when $\mathbf{X}$ is a path starting at $v_0$ or any
cyclic path. In the following example, we apply this test to finite
data sets simulated from a~lumped Markov chain.\vspace{-3pt}

\begin{exmp}[(Order estimation for a lumped reversible Markov chain)]
\label{example1}
A~random walk was simulated on the 9-state graph shown in Figure~\ref
{example}, from which we omitted self-edges on every state, all
weighted by 1. The observation was lumped into the 3 \textit{macrostates}
separated by the dashed lines. This is meant to illustrate a natural
experiment, where the difference between the states within each
macrostate is obscured by the measurement. From the resulting sequence,
we take the initial macrostate and every 7th macrostate thereafter to
form a path $\mathbf{X}$ of length 1000 in $\mathcal{X}=\{1,2,3\}$.

We test 4 reversible Markov models, that differ in the length of memory:
\begin{enumerate}
\item A first-order, reversible Markov chain.
\item A second-order, reversible Markov chain.
\item A variable-order model with maximum order 2, where states 1 and 3
are histories. Intuitively, only state 2 has ``memory.''
\item A variable-order model with maximum order 2, where states 2 and 3
are histories. Intuitively, only state 1 has ``memory.''
\end{enumerate}

For each model $i$, we assign a prior $\psi^{(i)}_{w,v_0}$ to the
transition matrix, where $v_0$ is the initial state in $\mathbf{X}$,
$w(u)=2$ for all $u\in\mathcal{X}^{3}$ and $\beta$ is the shortest
palindrome starting with $v_0$. We compared the 4 models using 50
independent realizations of the lumped Markov chain and found that
model 3 had the highest evidence in 72\% of the cases, while model 2
was selected in all the remaining cases. In Figure~\ref{boxplot}, we
report a boxplot of the logarithm of the Bayes factors comparing models
1, 2, and 4 against model 3.

%
\begin{figure}

\includegraphics{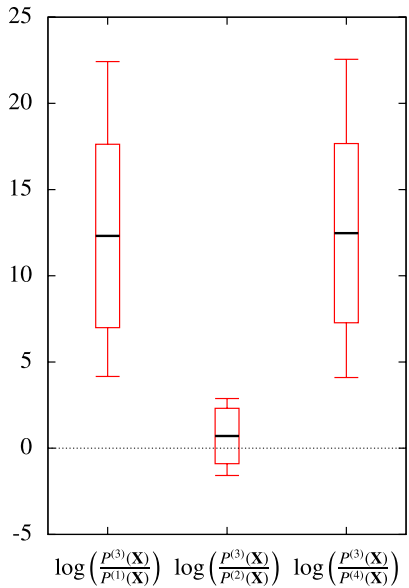}

\caption{Boxplot of logarithmic Bayes factors computed from 50
independent datasets.}
\label{boxplot}

\end{figure}

This represents compelling evidence for model 3. The result is not
entirely surprising given that this model gives memory to state 2,
which is slowly mixing, as indicated in Figure~\ref{example}. The fact
that the most complex model (model~2) is not necessarily selected
showcases the automatic penalty for model complexity in Bayes
factors.\vspace{-3pt}
\end{exmp}

%
\begin{figure}

\includegraphics{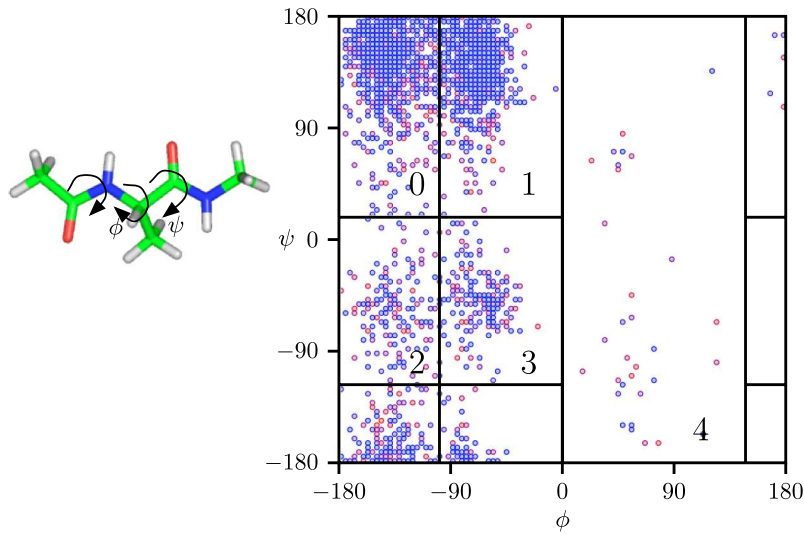}
\vspace{-3pt}
\caption{The structure of Ace-Ala-Nme is described by two dihedral
angles, $\phi$ and $\psi$. The periodic map on the right shows a
partition of conformational space into 5 states. The colored markers
indicate the free energy of bins centered at each point, which reveals
the metastable nature of this molecule's dynamics.}
\label{alanine}
\vspace{-5pt}
\end{figure}

We conclude this section with two applications of Bayesian analysis of
reversible Markov chains to molecular dynamics (MD). An MD simulation
approximates the time-reversible dynamics of a molecule in solvent. The
trajectories produced by a simulation are discretized in space and
time.\vspace{-3pt}

\begin{exmp}
The terminally blocked alanine dipeptide, shown in Figure~\ref
{alanine}, is a common test system for Markov models of MD. The
conformational space of the molecule, which is represented in the
figure in a two-dimensional projection, is partitioned into 5 states.
The states are believed to be metastable due to the basins that
characterize the free-energy function, also plotted in the figure. This
metastability allows one to approximate the dynamics of the molecule,
projected onto the partition, as a reversible Markov chain. The
approximation will be good when the discrete time interval at which a
trajectory is sampled is larger than the timescale for equilibration
within every state, but smaller than the timescale of transitions.

Few statistical validation methods are available for Markov models of
MD. Bacallado,
 Chodera and
 Pande  used a Bayesian hypothesis test to compare
different partitions of conformational space \cite{Bacallado2009ri}.
Here, we apply Bayes factors to test a first-order Markov model on a
fixed partition, by comparing it to second-order and variable-order
models on the same partition. The data $\mathbf{X}$ are the transition
counts in a single MD trajectory of 1767 steps sampled at an interval
of 6 picoseconds, as recorded in Table \ref{counts}. The prior
parameters $w$ and $\beta$ are the same as in the previous example. The
results of the model comparison are summarized in the following table.

%
\begin{table}
\tabcolsep=0pt
\vspace{-3pt}
\caption{Molecular dynamics simulation of the alanine dipeptide. The
entries in the table are the transition counts $(x_1,x_2)\to(x_2,x_3)$
in the trajectory $\mathbf{X}$, which has initial state $(0,4)$}
\label{counts}
\begin{tabular*}{\textwidth}{@{\extracolsep{\fill}}lcd{3.0}d{3.0}d{2.0}d{2.0}ccccd{2.0}d{2.0}cc@{}}
\hline
  &    & \multicolumn{5}{c}{$\bolds{x_3}$} &    &    &
\multicolumn{5}{c@{}}{$\bolds{x_3}$}\\[-5pt]
   &     & \multicolumn{5}{c}{\hrulefill}
 &     &    &
\multicolumn{5}{c@{}}{\hrulefill}\\
 $\bolds{x_1}$&$\bolds{x_2}$ & \multicolumn{1}{c}{\textbf{0}} & \multicolumn{1}{c}{\textbf{1}}& \multicolumn{1}{c}{\textbf{2}}
& \multicolumn{1}{c}{\textbf{3}}
& \multicolumn{1}{c}{\textbf{4}} & $\bolds{x_1}$& $\bolds{x_2}$& \multicolumn{1}{c}{\textbf{0}}
& \multicolumn{1}{c}{\textbf{1}}& \multicolumn{1}{c}{\textbf{2}}& \multicolumn{1}{c}{\textbf{3}}
& \multicolumn{1}{c@{}}{\textbf{4}}\\
\hline
 {0} & 0 & 261 & 187 & 13 & 2 & 0 &  3   &
0 &5 & 13& 2& 0& 0\\
& 1 & 188 & 144 & 13 & 11 & 0 &&1&5 & 4& 2& 1& 0\\
& 2 & 12 & 4 & 9 & 15& 0 &&2&4 & 3& 16& 5& 0\\
& 3 & 5& 1& 0& 1& 0 &&3&2 & 5& 3& 3& 0\\
& 4 & 1 & 0& 0& 0& 0 &&4&0 & 0& 0& 0& 0\\
[6pt]
 {1} & 0 &180 & 143& 22& 5& 0 &  4  & 0 &1
& 0& 0& 0& 0\\
&1&141& 125& 5& 5& 0 & &1&0 & 0 & 0& 0& 0 \\
&2&4 & 3& 10& 4& 0 & &2&0 & 0 & 0& 0& 0\\
&3&4 & 1& 10& 3& 0 & &3&0 & 0 & 0& 0& 0\\
&4&0 & 0& 0& 0& 0 & &4&0 & 0 & 0& 0& 0 \\
[6pt]
 {2} & 0 &16 & 13& 3& 0& 0 &\\
&1&16 & 4& 1& 1& 0 & \\
&2&12 & 12 & 37 & 11 & 0 &\\
&3&9 & 5& 15& 6& 0 &\\
&4&0 & 0& 0& 0& 0 & \\
\hline
\end{tabular*}
\vspace{-8pt}
\end{table}

\begin{center}
\vspace*{12pt}
\begin{tabular}{@{}lc@{}}
\hline
\textbf{Model ($\bolds i$)} & $\bolds{\log P^{(i)}(\mathbf{X})}$ \\
\hline
First order &  $-1846$ \\
Variable order 0 &  $-1824$ \\
Variable order 1 &  $-1825$ \\
Variable order 2 &  $-1844$ \\
Variable order 3 &  $-1846$ \\
Variable order 4 &  $-1847$ \\
Second order & \textbf{$\bolds-$1800} \\
\hline
\end{tabular}\vspace*{12pt}
\end{center}

The state describing each variable order model is the only state in the
model that has a memory of length 2 (the only state that is not a
history). There seems to be substantial evidence in favor of a
second-order model. Adding memory to states seen in a large number of
transition makes a bigger difference, as expected. This result is in
accordance with certain exploratory observations which indicate that at
the timescale of 6 picoseconds, the effect of water around the
molecule, neglected in our state definitions, persists.
\end{exmp}

\begin{exmp}
The alanine pentapeptide is a longer polymer that exhibits a higher
degree of structural and dynamical complexity. Buchete and Hummer
partition the conformational space of the molecule into 32 states by
chemical conventions \cite{Buchete2008}. An MD trajectory\footnote
{Simulated with the Amber-GSs forcefield at 300K in explicit solvent.}
in conformational space was projected onto this partition, and an
exploratory analysis suggested that the effects of memory decay after
500 picoseconds. Accordingly, we take a~conformation from the
trajectory every 500 picoseconds to form a sequence $\mathbf{X}$ of
1885 steps in $\mathcal{X}= \{0,\ldots ,31\}$.

As in previous examples, we tested models with varying lengths of
memo\-ry. Each model was assigned a conjugate prior, this time setting
$w(u) = 1/32$ for all $u\in\mathcal{X}^3$. Of all the variable-order
models where a single state has a memory of length 2 and all others are
histories, we found that only 4 models where strongly selected over a
first-order model. In the following table, we show the logarithm of the
evidence for each of these models, a~first-order model  and a
variable-order model that gives a memory of length 2 to all 4
states.

\begin{center}\vspace*{12pt}
\begin{tabular}{@{}lc@{}}
\hline
\textbf{Model ($\bolds i$)} & $\bolds{\log P^{(i)}(\mathbf{X})}$ \\
\hline
First order & $-4090.0$ \\
Variable order 14 & $-4015.5$ \\
Variable order 15 & $-3814.5$\\
Variable order 30 & $-3860.3$\\
Variable order 31 & $-3301.6$ \\
Variable order 14, 15, 30, 31 & \textbf{$\bolds-$2964.3} \\
\hline
\end{tabular}\vspace*{12pt}
\end{center}

This represents compelling evidence for a model that gives memory to
states 14, 15, 30  and 31. It is interesting to contrast inferences
based on this model to those based on a first-order Markov model. To do
this, we computed 1000 approximate posterior samples of the transition
matrix in each case. This was done by simulating a reinforced random
walk, which is a mixture of variable-order Markov chains with the
posterior distribution of $T$ as a mixing measure (see Proposition \ref
{posterior_sampling}). The reinforced random walk was simulated $10^7$
steps to obtain each sample.

%
\begin{figure}

\includegraphics{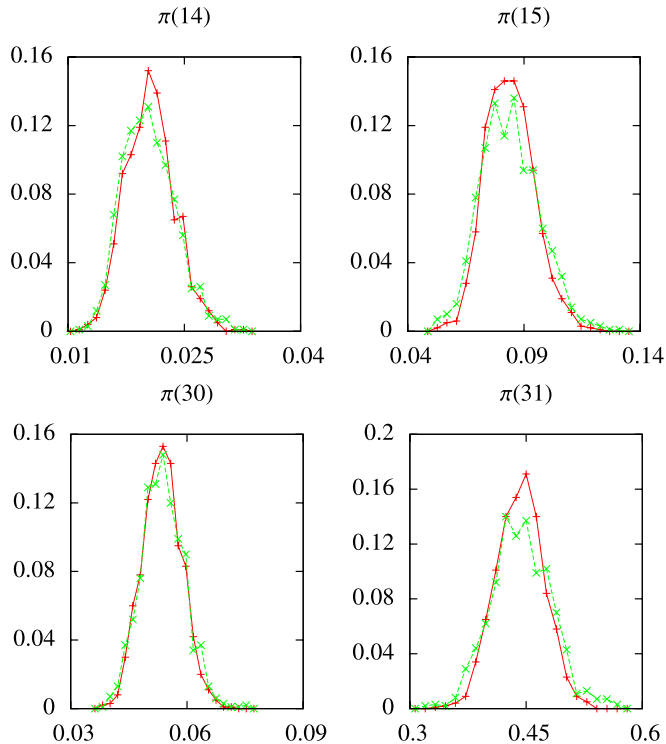}

\caption{Histograms of 1000 posterior samples of the stationary
probabilities of states 14, 15, 30, 31. The red solid lines correspond
to the first-order Markov model, and the green dashed lines to the
variable-order Markov model that gives a memory of length 2 to states
14, 15, 30  and 31.}
\label{histograms_stat}

\end{figure}
%
\begin{figure}[t]
\vspace{-3pt}
\includegraphics{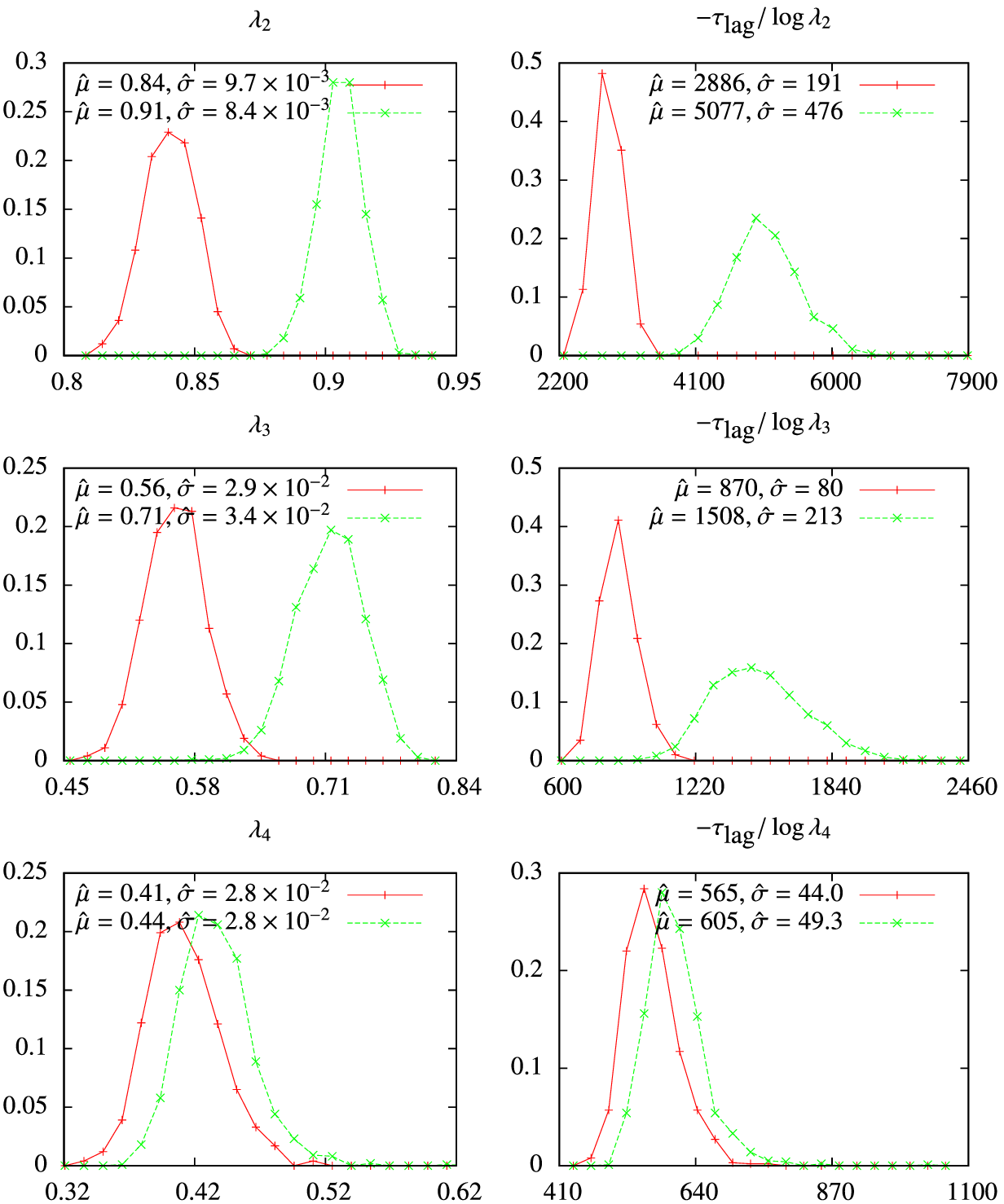}
\vspace{-6pt}
\caption{Histograms of 1000 posterior samples of the second, third and
fourth largest eigenvalues of the transition matrix, as well as the
timescales associated with these eigenvalues. The red solid lines
correspond to the first-order Markov model, and the green dashed lines
to the variable-order Markov model that gives memory to states 14, 15,
30  and 31. In both cases, we compute the eigenvalues of the transition
matrix for the process $V_n, n\in\mathbb{N}$, in $\mathcal{X}^2$. All
sample means $\hat{\mu}$ and standard deviations $\hat{\sigma}$ are shown.}
\label{histograms eig tim}%
\vspace{-10pt}
\end{figure}

In Figure~\ref{histograms_stat}, we histogram stationary probabilities
of the transition matrices sampled from the posterior. In particular,
we show plots for the stationary probabilities of states 14, 15, 30
and 31. In the variable-order model, we define $\pi(x) = \sum_{y\in
\mathcal{X}} \pi(xy)$. The inferences of each model in this case are
very similar.

The largest eigenvalues of the transition matrix are also of interest
because they are related to different modes of relaxation. Each
eigenvalue $\lambda$ is associated with a timescale $-\tau_
{\mathrm{lag}}/\log\lambda$, which is useful in exploratory analysis. Here, $\tau
_{\mathrm{lag}}$ is the length in time of one step of the Markov chain, or
500 picoseconds. In Figure~\ref{histograms eig tim}, we histogram
posterior samples of the three largest nonunit eigenvalues and their
associated timescales. In this case, the inferences of each model are
quite different, with the variable-order model predicting larger
eigenvalues and timescales.
\end{exmp}

\section{Conclusions}\label{sec6}
We define a reinforcement scheme for the higher-order, reversible
Markov chain that extends the ERRW on an undirected graph. Several
properties of the ERRW, like recurrence and partial exchangeability,
were shown to generalize to this process. Other properties may also
generalize but were not pursued here. In particular, we can mention the
uniqueness results of   Johnson \cite{Zabell1982ir} and Rolles
\cite{Rolles2003xd}, and the fact that mixtures of measures in
$\mathcal{D}$ are weak-star dense in the space of all priors \cite
{Diaconis2006gd}.

The reinforced random walk leads to a conjugate prior that facilitates
estimation and hypothesis testing of reversible processes in which the
effects of memory decay after some time. Certain statistical problems
remain a challenge, such as inferring the transition matrix with a \textit
{fixed} stationary distribution. In applications, it will become
important to evaluate the objectivity of the prior and to determine the
optimal value of its parameters in this sense.

From a practical point of view, we only discussed Bayesian updating for
data sets composed of a single Markov chain starting with probability 1
from the initial state $v_0$ used in the prior. Numerical algorithms
are needed to perform inference with data sets composed of multiple
chains. A starting point could be the method developed by Bacallado,
 Chodera and
 Pande  to apply the prior of Diaconis and Rolles to first-order,
reversible Markov chains \cite{Bacallado2009ri}.

\begin{appendix}
\section*{Appendix}\label{appm}

In the following, we use the notation defined in the first paragraph of
Section \ref{sec2}.

\begin{prop}[(Kolmogorov's criterion)]
Let $X_n, n\in\mathbb{N}$, be an irreducible order-$r$ Markov chain with
transition probabilities $p$. Then $X_n$ is reversible if and only if
for any cyclic admissible path $v_0,v_1,\ldots ,v_n,v_0$,
%
\begin{equation}
\label{kolmogorov}
p(v_1|v_0)p(v_2|v_1)\cdots p(v_0|v_n) =
p(v_0^*|v_1^*)p(v_1^*|v_2^*)\cdots p(v_n^*|v_0^*).
\end{equation}
\end{prop}

\begin{pf}
The ``only if'' statement is straightforward. By the definition of the
stationary distribution and reversibility
\begin{eqnarray*}
p(v_1|v_0)p(v_2|v_1)\cdots p(v_0|v_n) &=& \frac{P_\pi(\overline
{v_0v_1})}{\pi(v_0)}\frac{P_\pi(\overline{v_1v_2})}{\pi(v_1)}\cdots\frac
{P_\pi(\overline{v_nv_0})}{\pi(v_n)} \\
&=& \frac{P_\pi(\overline{v_1^*v_0^*})}{\pi(v_0^*)}\frac{P_\pi(\overline
{v_2^*v_1^*})}{\pi(v_1^*)}\cdots\frac{P_\pi(\overline{v_0^*v_n^*})}{\pi
(v_n^*)} \\
&=& p(v_0^*|v_1^*)p(v_1^*|v_2^*)\cdots p(v_n^*|v_0^*).
\end{eqnarray*}
To prove the ``if'' statement, choose an arbitrary state $u$; then, for
any $v$, since the chain is irreducible, there is an admissible path
$u,v_1,v_2,\ldots ,v_n,v$ with positive probability. Define
%
\begin{equation}
\label{stat_construction}
\pi'(v) = B\frac{p(v_1|u)p(v_2|v_1)\cdots
p(v|v_n)}{p(v_n^*|v^*)p(v_{n-1}^*|v_n^*)\cdots p(u^*|v_1^*)},
\end{equation}
where $B$ is a positive constant. Note that this expression does not
depend on the sequence $v_1,\ldots ,v_n$ chosen. Take a different
sequence $z_1,\ldots ,z_m$. Let $t\in\mathcal{X}^r$ be a palindrome,
then because the chain is irreducible, we can find a sequence
$v,t_1,t_2,\ldots ,t$ with positive probability, and it is easy to see
from equation~\eqref{kolmogorov} that the palindrome $v,t_1,t_2,\ldots
,t,\ldots ,t_2^*,t_1^*,v^*$ has positive probability. We can construct
another palindrome $u^*,s_1,s_2,\ldots ,s_2^*,s_1^*,u$ in the same way.
Multiplying equation~\eqref{stat_construction} by factors of 1,
\begin{eqnarray*}
B\frac{p(v_1,v_2,\ldots ,v|u)}{p(v_n^*,v_{n-1}^*,\ldots ,u^*|v^*)}
&=  & B \frac{p(v_1,v_2,\ldots ,v|u)}{p(v_n^*,v_{n-1}^*,\ldots ,u^*|v^*)}
\frac{p(t_1,t_2,\ldots ,v^*|v)}{p(t_1,t_2,\ldots ,v^*|v)} \\
&&{} \times \frac{p(z_m^*,z_{m-1}^*,\ldots ,u^*|v^*)}{p(z_1,z_2,\ldots
,v|u)}\frac{p(s_1,s_2,\ldots ,u|u^*)}{p(s_1,s_2,\ldots ,u|u^*)} \\
&&{}\times\frac{p(z_1,z_2,\ldots ,v|u)}{p(z_m^*,z_{m-1}^*,\ldots
,u^*|v^*)} \\
&=  & B\frac{p(z_1,z_2,\ldots ,v|u)}{p(z_m^*,z_{m-1}^*,\ldots ,u^*|v^*)}.
\end{eqnarray*}
The first four terms equal 1 because the numerator and denominator are
the probabilities of the same cycle forward and backward, which are
equal by equation~\eqref{kolmogorov}. Now, we check that $\pi'(v)$ satisfies
the reversibility conditions specified in the \hyperref[sec1]{Introduction}. First, we
show that $\pi'(v)=\pi'(v^*)$. Take a path $u,z_1,\ldots ,z_\ell,v^*$
with positive probability, and the previously found palindrome
$u^*,s_1,s_2,\ldots ,s_2^*,s_1^*,u$, then applying the same method,
\begin{eqnarray*}
\pi'(v) &=& B\frac{p(v_1,v_2,\ldots ,v|u)}{p(v_n^*,v_{n-1}^*,\ldots
,u^*|v^*)}\\
&=&   B \frac{p(v_1,v_2,\ldots ,v|u)}{p(v_n^*,v_{n-1}^*,\ldots ,u^*|v^*)}
\frac{p(s_1,s_2,\ldots ,u|u^*)}{p(s_1,s_2,\ldots ,u|u^*)} \\
&&{}\times\frac{p(z_\ell^*,z_{\ell-1}^*,\ldots ,u^*|v)}{p(z_1,z_2,\ldots
,v^*|u)} \frac{p(z_1,z_2,\ldots ,v^*|u)}{p(z_\ell^*,z_{\ell-1}^*,\ldots
,u^*|v)} \\
&=  & B\frac{p(z_1,z_2,\ldots ,v^*|u)}{p(z_\ell^*,z_{\ell-1}^*,\ldots
,u^*|v)} =\pi'(v^*).
\end{eqnarray*}
From this, and equation~\eqref{stat_construction} we deduce that for any
admissible $v,z$, $\pi'(v)p(z|\break v) = \pi'(z^*)p(v^*|z^*)$. Since the
state space is finite, we can choose $B$ such that $\pi'$ sums to 1. We
have shown that the weights $k_{v,z} \equiv\pi'(v)p(z|v)$ satisfy the
conditions of a reversible random walk with memory, so by Proposition
\ref{representation} the process with transition probabilities $p$
represents a reversible, order-$r$ Markov chain.
\end{pf}

\begin{pf*}{Proof of Proposition \ref{mixed-PE}}
The probability $H_{w,v_0}(v_0,\ldots ,v_m)$ is a~pro\-duct of transition
probabilities, to which the $n$th transition contributes a~factor of
%
\begin{equation}
p_n = \frac{w'_{n-1}(\overline{f(v_{n-1})v_n})}{w''_{n-1}(f(v_{n-1}))}.
\label{factor1}
\end{equation}
We know that $f(v_n)$ cannot be longer than $\overline{f(v_{n-1})v_n}$
by Proposition \ref{regaining_memory}; let $L(v_{n-1},v_n)$ be the set
of histories of $v_n$ that are shorter than $f(v_{n-1})$. If this set
is nonempty, let us multiply equation~\eqref{factor1} by factors of 1, to
obtain the following factor for the $n$th transition:
%
\begin{equation}
p_n = \frac{w'_{n-1}(\overline{f(v_{n-1})v_n})}{w''_{n-1}(f(v_{n-1}))}
\prod_{z\in L(v_{n-1},v_n)}\frac{w'_{n-1}(z^+)}{w''_{n-1}(z^+)},
\label{factor2}
\end{equation}
where $z^+$ is the ending of $v_n$ that is longer than $z$ by 1. The
added factor equals 1 because, if $w'_{n-1}(z^+)\neq w''_{n-1}(z^+)$,
then $f(v_{n-1})$ must end on $z$, which by definition is a history
shorter than $f(v_{n-1})$, a contradiction.

Consider all the possible factors in the numerator of
$H_{w,v_0}(v_0,\ldots ,v_m)$. Take any $h\in\mathscr{H}$ that is \textit
{minimal}, meaning  that it  does not end in another history. For any
$a\in\mathcal{X}$, we will see a factor $w'(ha)$ after every transition
through $ha$. The conjugate factor $w'(ah^*)$ will appear every time we
go through $ah^*$, because:
\begin{itemize}
\item
If $\mathrm{A}(ah^*)\in\mathscr{H}$, it is minimal by the closure
properties of $\mathscr{H}$, so $w'(ah^*)$ will be the numerator of the
first factor in equation~\eqref{factor2}.
\item
Otherwise, the minimal history in the transition ending in $ah^*$ will
be longer than $h^*$, and there will be an added factor in equation~\eqref
{factor2} with $w'(ah^*)$ in the numerator. Conversely, note that the
factor $w'(ah^*)$ is only added to the numerator of equation~\eqref{factor2}
when we go through $ah^*$ for some minimal $h$, because we required
that $h^*\in L(v_{n-1},v_n)$, so $h\in\mathscr{H}$ and does not end in
another history.
\end{itemize}
As in the proof of Proposition \ref{PE}, we argue that every new factor
$w'(ha)$ or $w'(ah^*)$ is increased by $c$ with respect to the previous
one (or by $2c$ if $ha$ is a palindrome). Therefore, the numerator of
$H_{w,v_0}(v_0,\ldots ,v_m)$ is only a function of the transition counts
and the initial state.

Finally, consider all the factors in the denominator of
$H_{w,v_0}(v_0,\ldots ,v_m)$. Take any minimal history $h$. We will see
a factor $w''(h)$, for every transition through $h$. The conjugate
factor $w''(h^*)$ will appear every time we go through $h^*$, because:
\begin{itemize}
\item
If $h^*$ is also minimal, then $w''(h^*)$ will be in the denominator of
the first factor in equation~\eqref{factor2}.
\item
Otherwise, we know that $\mathrm{A}(h^*)$ is not a history, so the
transition ending in $h^*$ must have a history at least as long as
$h^*$, which is longer than the history $\Omega(h^*)$. So, $w''(h^*)$
will appear in the denominator of a factor added in equation~\eqref{factor2}.
Conversely, we only add factors of $w''(h^*)$ to the denominator of
equation~\eqref{factor2} when we go through $h^*$ for a minimal $h$, because
we required $\Omega(h^*)\in L(v_{n-1},v_n)$ which implies $h$ minimal.
\end{itemize}
As before, every new factor $w''(h)$ or $w''(h^*)$ will be increased by
$c$ with respect to the previous one (or by $2c$ if $h$ is a
palindrome). Therefore, the denominator is a function of the transition
counts and the initial state, and the process is partially exchangeable.
\end{pf*}

\begin{pf*}{Proof of Proposition \ref{Concentration}}
Let $\vec{C}_n(u,v)$ be the transition counts from $u$ to $v$ in the
first $n$ steps of a stochastic process on $\mathcal{X}^r$. Also,
define $\vec{C}_n(u)\equiv\sum_{v\in\mathcal{X}^r}\vec{C}_n(u,v)$,
which counts the visits to $u$. Remember $\mathcal{T}''$ is the set of
irreducible transition matrices for variable-order, reversible Markov
chains where all $h\in\mathscr{H}$ are histories. Define the event $D$,
that the set $\{\vec{C}_n(u,v)/\vec{C}_n(u)\dvtx \forall u,v
$ admissible$\}$ converges to a transition probability matrix in $\mathcal{T}''$.

From the recurrence of the variable-order, reinforced random walk and
equation~\eqref{mixed-order mixture}, it is evident that the set of
irreducible Markov chains has measure 1 under $\psi_{w,v_0}$. In this
set, the variables $\{\vec{C}_n(u,v)/\vec{C}_n(u)\dvtx \forall u,v \mbox{
admissible}\}$ converge almost surely to the transition probabilities,
so for any $T\notin\mathcal{T}''$ irreducible, \mbox{$P_{v_0}^T(D)=0$}.
Furthermore, by Lemma \ref{Q}, $D$ happens almost surely in the
variable-order, reinforced random walk. Putting this into equation~\eqref
{mixed-order mixture}, we have
\begin{eqnarray*}
H_{w,v_0}(D) = 1 = \int_\mathcal{T} P_{v_0}^T(D)\,d\psi_{w,v_0}(T) \leq
\int_{\mathcal{T}''}\,d\psi_{w,v_0}(T),
\end{eqnarray*}
which implies the proposition.
\end{pf*}

\begin{lem}[(L\'{e}vy)] \label{levy} Consider a sequence of events $B_k\in
\mathcal{F}_k, k\in\mathbb{N}$, in some filtration $\{\mathcal{F}_k\}$.
Let $b_n = \sum_{k=1}^n \mathbf{1}_{B_n}$ be the total number of events
occurring among the first $n$, and let $s_n = \sum_{k=1}^n
P(B_k|\mathcal{F}_{k-1})$ be the sum of the first $n$ conditional
probabilities. Then, for almost every $\omega$:
\begin{itemize}
\item If $s_n(\omega)$ converges as $n\to\infty$, then $b_n(\omega)$
has a finite limit.
\item If $s_n(\omega)$ diverges, then $b_n(\omega)/s_n(\omega)\to1$.
\end{itemize}
\end{lem}

\begin{lem}
$H_{w,v_0}(D) = 1$.
\label{Q}
\end{lem}

\begin{pf}
For any $u$ in $\{\mathcal{X}^q\dvtx q\leq r+1\}$, the variables
$n^{-1}w'_n(u)$ and\break $n^{-1}w''_n(u)$ are functions of $\{n^{-1}\vec
{C}_n(u,v)\dvtx \forall u,v \mbox{ admissible}\}$, therefore they converge
almost surely, because the reinforced random walk is a mixture of
irreducible Markov chains for which the latter converge. The
reinforcement scheme defined in Definition~\ref{mixed-order reinforcement} imposes
some constraints on the limits of $n^{-1}w'_n(u)$ and $n^{-1}w''_n(u)$.
Note that $w''_n(u)$, $w''_n(u^*)$, $w'_n(u)$  and $w'_n(u^*)$ never
differ by more than $c$; we also know that the reinforced random walk
is positive recurrent (it is a mixture of irreducible, finitely-valued
Markov chains), so almost surely
%
\begin{eqnarray}
\label{reversibility of limit}
\lim_{n\to\infty}n^{-1}w''_n(u) &= &\lim_{n\to\infty} n^{-1}w''_n(u^*)=\lim_{n\to\infty}n^{-1}w'_n(u)
\nonumber
\\[-8pt]
\\[-8pt]
 &=& \lim_{n\to\infty} n^{-1}w'_n(u^*)>0.
\nonumber
\end{eqnarray}
Denote this limit $w_\infty(u)=w_\infty(u^*)$. It is also easy to see
that if $u\in\mathcal{X}^q$, then for all $s>q$,
%
\begin{equation}
\label{balance}
\sum_{\{v\in\mathcal{X}^s\dvtx  v \mathrm{\ ends\ in\ } u\}} w_\infty(v) =
w_\infty(u).
\end{equation}

Now, let $\tau_n$ be the $n$th visit to $u\in\mathcal{X}^r$ and let
$B_n$ be the event that we make a transition to $v$ at $\tau_n$. Define
\[
p_n(f(u),v) \equiv H_{w,v_0}(B_n|\sigma(Y_1,\ldots ,Y_{\tau_n}))= \frac
{w'_{\tau_n}(\overline{f(u)v})}{w''_{\tau_n}(f(u))}.
\]
We know $p_n(f(u),v)$ converges a.s.  to $w_\infty(\overline
{f(u)v})/w_\infty(f(u))>0$. Therefore, $\sum_n p_n(f(u),v) = \infty$
a.s., and by L\'{e}vy's extension of the Borel--Cantelli lemma (Lemma
\ref{levy}),
\begin{eqnarray*}
&&\frac{\sum_{m=1}^n \mathbf{1}_{B_m}}{\sum_{m=1}^n p_m(f(u),v)}  \to1
  \quad \mbox{a.s.} \\
 &&\qquad \implies \quad \lim_{n\to\infty}\frac{1}{n} \sum_{m=1}^n \mathbf{1}_{B_m} =
\lim_{k\to\infty} \frac{\vec{C}_k(u,v)}{\vec{C}_k(u)} = \frac{w_\infty
(\overline{f(u)v})}{w_\infty(f(u))}.
\end{eqnarray*}
This means that $\{\vec{C}_n(u,v)/\vec{C}_n(u)\dvtx \forall u,v \mbox{
admissible}\}$ converges $H_{w,v_0}$-a.s.  to a set of transition
probabilities, $w_\infty(\overline{f(u)v})/w_\infty(f(u))$, for a
variable-order Markov chain with histories $\mathscr{H}$. To show that
this Markov chain is reversible, note that $w_\infty$ is the stationary
distribution, because
\begin{eqnarray*}
\sum_{
{\fontsize{8.36pt}{10pt}\selectfont{{\left \{\matrix{u\in\mathcal{X}^r{:} \cr u,v \mathrm{\
admissible}} \right\}}}}} w_\infty(u)\frac{w_\infty(\overline
{f(u)v})}{w_\infty(f(u))} & =&
\sum_{ {\fontsize{8.36pt}{10pt}\selectfont{ \left \{\matrix{h\in\mathscr{H} \mathrm{\ minimal{:}} \cr h,v
\mathrm{\
admissible}} \right\}}}} \sum_{\{u\dvtx f(u)=h\}} w_\infty(u) \frac{w_\infty
(\overline{hv})}{w_\infty(h)} \\
& =& \sum_{ {\fontsize{8.36pt}{10pt}\selectfont{\left \{\matrix{h\in\mathscr{H} \mathrm{\ minimal{:}} \cr h,v
\mathrm{\ admissible}} \right\}}}} w_\infty(h) \frac{w_\infty(\overline
{hv})}{w_\infty(h)} \\
& =& w_\infty(v),
\end{eqnarray*}
where we used equation~\eqref{balance} in the last two identities. By equation~\eqref{reversibility of limit}, $w_\infty$ satisfies the conditions for
reversibility. Therefore, $H_{w,v_0}(D)=1$.
\end{pf}
\end{appendix}

\section*{Acknowledgments}
The author would like to thank Persi Diaconis and Vijay Pande for
valuable suggestions, and Lutz Maibaum for providing molecular dynamics
datasets. 

\begin{supplement}[id=supp]
\stitle{Law of a variable-order, reinforced random walk\\}
\slink[doi,text={10.1214/10-AOS857 SUPP}]{10.1214/10-AOS857SUPP} 
\sdatatype{.pdf}
\sfilename{supplement.pdf}
\sdescription{We provide a closed form expression for this law as a
function of
transition counts and suggest how it could be useful.}
\end{supplement}

\vfill\eject
\printaddresses

\end{document}